\documentclass[11pt,twoside]{article}
\usepackage{latexsym}
\usepackage{amssymb,amsbsy,amsmath,amsfonts,amssymb,amscd}
\usepackage{graphicx,color}
\usepackage{hyperref}
\newcommand{\commentout}[1]{}
\newcommand{\R}{\mathbb{R}}

\newcommand {\e}  {\varepsilon}

\newcommand {\lb} {\lambda}
\newcommand {\Chi} {{\bf \raise 2pt \hbox{$\chi$}} }
\newcommand {\f}   {\frac}
\newcommand {\p}   {\partial}
\newcommand{\fer}{\eqref}
\newcommand{\dis}{\displaystyle}
\newcommand {\proof} {\noindent {\bf Proof}. }
\newcommand{\beq}{\begin{equation}}
\newcommand{\eeq}{\end{equation}}
\newcommand{\bea} {\begin{array}{rl}}
\newcommand{\eea} {\end{array}}
\newcommand{\bepa}{\left\{ \begin{array}{l}}
\newcommand{\eepa} {\end{array}\right.}
\newtheorem{theorem}{Theorem}[section]
\newtheorem{lemma}[theorem]{Lemma}

\newcommand{\qed}{{ \hfill
                       {\unskip\kern 6pt\penalty 500 \raise -2pt\hbox{\vrule\vbox to 6pt{\hrule width 6pt
                       \vfill\hrule}\vrule} \par}   }}
\title{Asymptotic analysis of a selection model with space}

\author{
Sepideh Mirrahimi\thanks{Institut de Math\'ematiques de Toulouse; UMR 5219, Universit\'e de Toulouse; CNRS, UPS IMT, F-31062 Toulouse Cedex 9, France; E-mail:Sepideh.Mirrahimi@math.univ-toulouse.fr} 
\and Beno\^ \i t Perthame\thanks{Sorbonne Universit\'es, UPMC Univ. Paris 06, UMR 7598, Laboratoire Jacques-Louis Lions, F-75005, Paris, France}
\thanks{CNRS, UMR 7598, Laboratoire Jacques-Louis Lions, F-75005, Paris, France
} 
\thanks{INRIA Paris-Rocquencourt, EPC Mamba, Institut Universitaire de France. Email: benoit.perthame@upmc.fr} 
}

\date{\today}

\begin{document}
\maketitle
\pagestyle{plain}
%\tableofcontents
\pagenumbering{arabic}

\begin{abstract}
Selection of a phenotypical trait can be described in mathematical terms by `stage structured' equations which are usually written under the form of integral equations so as  to express competition for resource between individuals whatever is their trait. The solutions exhibit a concentration effect  (selection of the fittest); when a small parameter is introduced they converge to a Dirac mass. 

An additional  space variable can be considered in order to take into account local environmental conditions. Here we assume this environment is a single nutrient which diffuses in the domain. In this framework, we prove that the solution converges to a Dirac mass in the physiological trait which depends on time and on the location in space with Lipschitz continuity. The main difficulties come from the lack of compactness in time and trait variables. Strong convergence can be recovered from uniqueness in the limiting constrained equation after Hopf-Cole change of unknown.

Our analysis is motivated by a model of tumor growth introduced in \cite{AL.TL.JC.AE.BP:14} in order to explain emergence of resistance to therapy. 

\end{abstract}

\bigskip

\\[5pt]

%%%%%%%%%%%%%%%%%%%%%%%%%%%%%%%%%%%%%%%%%%%%
\section{Setting the problem}
\label{sec:pb}
%-------------------------------------------
%%%%%%%%%%%%%%%%%%%%%%%%%%%%%%%%%%%%%%%%%%%%
%In this article, we study a model 

In this paper, we are interested in the study of the evolutionary dynamics of populations structured by phenotypical traits and space. While our motivation comes from the study of tumor growth, we investigate the selection of the fittest individuals and the heterogeneity of the population. A population structured by a phenotypical trait can be modeled using integro-differential Lotka-Volterra equations. The solutions of such equations, when we consider small mutation steps  and in long time, converge to Dirac masses (see \cite{GB.BP:08,AL.SM.BP:10}); this property corresponds to  the selection of the fittest traits. In this paper, we study such behavior considering a spatial structure for the population.
\\

\noindent
A simple way to describe the selection of the fittest  individuals, when environmental conditions depend on space, was proposed in \cite{AL.TL.JC.AE.BP:14} as a model for emergence of resistance to drug in cancer therapy. This model assumes the  evolution of cells and is written as a coupled system of integro-differential equations structured by trait~$x$ and by a space variable~$y$
\begin{equation} \begin{cases}
\e \partial_t n_\e(y,x,t) = \left[ r(x)  c_\e(t,y)  - d(x) (1+\varrho_\e(y,t)) \right] n_\e(y,x,t),  \qquad y \in \R, \; 0< x <1 , \; t \geq 0,
\\[8pt]
-\Delta_y c_\e (y,t) + [\varrho_\e(y,t) +\lb ] \; c_\e(y,t) = \lb c_B ,
\qquad 
\varrho_\e(y,t) = \int n_\e (y,x,t) dx .
\end{cases}
\label{eq:ntrue}\end{equation}
 The first equation describes the dynamics of a cell population density $n_\e$. The second equation describes a nutrient $c_\e$ (and a drug can be included in the same way) diffused within the tumor from a constant  input concentration $c_B$ with rate $\lb$.  
The term $r(x)$ denotes the proliferation rate of cells expressing trait $x$ due to the consumption of resource. 
The function $d(x)$ models the death rate of cells with trait $x$ due to the competition with other cells at the same position. The small parameter $\e$ is introduced to consider the long time behavior of the cell population. Note that we do not consider mutations in this model, supposing that all traits are already present in the population, possibly at very small quantities.
 \\

\noindent
Our goal is to show that, when $\e$ vanishes, there is selection of a space and time dependent fittest trait $X(y,t)$ in the cell population as numerically shown in \cite{AL.TL.JC.AE.BP:14}. 
\\

\noindent
In order to get more complete results, and  show better the difficulties when handling the time variable, we also  study a related model where the integral equation for $n_\e$ is coupled to a parabolic equation for the nutrient.
\begin{equation}\label{eq:ne}
\e \partial_t n_\e(y,x,t) = \left[ r(x)  c_\e(t,y)  - d(x) (1+\varrho_\e(y,t)) \right] n_\e(y,x,t),  \qquad y \in \R, \; 0< x <1 , \; t \geq 0,
\end{equation}
\begin{equation}\label{eq:ce}
\f{\p}{\p t} c_\e-\Delta_y c_\e (y,t) + [\varrho_\e(y,t) +\lb ] \; c_\e(y,t) = \lb c_B ,
\end{equation}
\begin{equation}\label{eq:rho}
\varrho_\e(y,t) = \int n_\e (y,x,t) dx .
\end{equation}

\noindent
Recent technologic advances reveal evidence of heterogeneity within cancer tumors (see for instance \cite{MG.AR.SH.JL:12}).
Taking into account  this intratumor heterogeneity is crucial in the study of the tumor growth and the emergence of drug resistance  (see \cite{TB.XX.GY.EW:09,CS:10,MG.AR.SH.JL:12}), and leads to important challenges in finding effective treatment strategies. 
The above model introduces a simple  way to include spatial and phenotypical structure of the cell population together with the diffusion of the nutrient in the domain. Our study  indicates that intratumor heterogeneity can emerge as an evolutionary process and provides a description of the space-dependent dominant traits.\\

\noindent
The dynamics of phenotypically structured populations under the effect of mutations and competition between the traits has been studied widely during the last decade using stochastic methods and integro-differential equations (see for instance \cite{OD.PJ.SM.BP:05,NC.RF.SM:06,LD.PJ.SM.GR:08,Champagnat.phd,Raoulphd,Mirrahimi.phd} and the references therein). In particular an approach based on a WKB type ansatz, leading to Hamilton-Jacobi equations, provides an analysis of the asymptotic behavior of the populations structured by phenotypical traits (see \cite{GB.BP:08,GB.SM.BP:09,AL.SM.BP:10} and the references therein). It can be shown using this method that in long time and considering small mutation steps the population concentrates on a dominant trait that evolves in time. In other words, the population density tends to a Dirac mass in the phenotypical trait which depends on time.\\

\noindent
Several refinements are necessary to consider problems which are more relevant biologically.  In particular, one can include the interaction of the population with nutrients (see \cite{SM.BP.JW:10,NC.PJ.SM:14} for some results in this direction). The introduction of the resource leads to the study of integro-differential coupled systems. 
Moreover, most of the cited works neglect the spatial structure of the environment and consider a well-mixed population. However, as mentioned above the environmental heterogeneity is an important element to be considered. The study of models of populations  structured jointly by space and trait
has gained much attention recently and leads to several important difficulties to be  overcome. These difficulties are mainly due to the integral term in one of the variables in the equations.
Most of the recent attempts to tackle these problems  concentrate mostly on the spatial propagation of the population and less on the phenotypical selection (see \cite{SM.GR:11,MA.JC.GR:13, EB.SM:14,EB.VC:14}, and also \cite{NC.SM:07,AA.LD.CP:12} for the stochastic derivation of such models and the study of steady states). 
In this paper, using the WKB approach mentioned above, we study models which take into account the spatial and phenotypical structure of the population and lead to the selection of a space and time dependent phenotypical trait. 
Note that in our model the spatial heterogeneity is induced by  a nutrient which diffuses in the domain. 
\\

\noindent
We assume some conditions on the model parameters and on the initial data 
\begin{equation}\label{as:0}
\varrho_m \leq \varrho_\e^0 (y) \leq \varrho_M,  \qquad \varrho_m := c_B \f{\lb}{\lb +\varrho_M} \,  \dis \min_{0< x<1} \f {r(x)}{d(x)} -1   >0,  \quad  \varrho_M := c_B \max_{0< x<1} \f {r(x)}{d(x)} -1,
\end{equation}
\begin{equation}\label{as:1}
c_m < c^0 (y) < c_B , \quad c_m := c_B \f{\lb}{\lb +\varrho_M}, \qquad c^0 \in W^{1,\infty} (\R) .
\end{equation}
Note that  the non-extinction condition, $\varrho_m>0$,  is equivalent to write
$$ 
\lb c_B \min_{0< x<1} \f {r(x)}{d(x)} \geq \lb +c_B \max_{0< x<1} \f {r(x)}{d(x)},
$$
a condition which is satisfied for parameters $\lb$ and $c_B$ sufficiently large.
\\

\noindent
For $n_\e$, we assume initially a `Gaussian type' concentration 
%------------------------------
\begin{equation}\label{as:2}
\begin{cases}
n^0_\e \rightharpoonup \varrho^0(y) \delta (x-X^0(y)), \quad \text{with the condition }Ê\ r(X^0(y)) c^0(y) - d(X^0(y)) (1+\varrho^0(y))=0,
\\[8pt]
n^0_\e = e^{u^0_\e(y,x)/\e} , \qquad u^0_\e(y,x) \underset{ \e \to 0}{\longrightarrow} u^0(y,x) \ \text{ locally uniformly},
\\[8pt]
| \f{\p}{\p x} u_\e^0|  +| \f{\p^3}{\p x^3} u_\e^0|  \leq K^0,   \qquad \f{\p^2}{\p x^2} u_\e^0 \leq - a <0.
\end{cases}
\end{equation}
In particular these conditions imply that 
$$
\max_{0\leq x \leq 1} u^0_\e(y,x) = u^0_\e(y, X_\e^0(y))  \underset{ \e \to 0}{\longrightarrow} 0 = \max_{0\leq x \leq 1} u^0(y,x) = u^0(y, X^0(y)) 
$$
with $ X_\e^0(y)) \underset{ \e \to 0}{\longrightarrow} X^0(y)$ locally uniformly. 

\noindent
Finally, we assume that $c$ and $d$ are smooth and that for some constant $K^0$
\begin{equation}\label{as:3} 
\begin{cases}
|r'|+ |d'| + |r''|+|d''| +|r'''|+ |d'''|  \leq K^0, 
\\[8pt]
r'' <0, \qquad d'' >0 .
\end{cases}
\end{equation}

%-------------------------------
\begin{theorem} [Parabolic case]  With assumptions \eqref{as:0}--\eqref{as:3}, there is $X(y,t) \in W^{1,\infty}(\R\times \R^+)$,  $\varrho(y,t)\in C( \R\times \R^+) $ and $c(y,t) \in C( \R\times \R^+) $, such that solutions of \eqref{eq:ne}--\eqref{eq:rho} satisfy, as $\e \to 0$, 
\begin{eqnarray}
\varrho_\e \to \varrho(y,t), \ \text{ almost everywhere}, 
\\
c_\e \to  c(y,t), \ \text{ locally uniformly}, 
\\
 n_\e(y,x,t) \rightharpoonup  \varrho(y,t) \delta \big(x-X(y,t) \big) , \ \text{ weakly in measures}.
\end{eqnarray}
Moreover, we have 
\begin{equation}\label{eq:ctrue}
\f{\p}{\p t} c-\Delta_y c (y,t) +  [ \, \varrho (y,t) +\lb \, ] \; c(y,t) = \lb c_B , \qquad y \in \R, \; t\geq 0 ,
\end{equation}
and
\begin{equation}\label{eq:X}
 r\big(X(y,t)\big) c(y,t) - d\big(X(y,t)\big) \big(1+\varrho(y,t)\big)=0.
 \end{equation}
Finally, $\varrho(y,t)$  and $c(y,t) $ are H\"older continuous in $t$ and Lipschitz continuous in $y$.
\label{th:1}
\end{theorem}
%-------------------------------

\noindent
The convergence of $c_\e$ can be derived from parabolic regularity, while more elaborate arguments are needed to obtain the limit of $n_\e$. To obtain a priori bounds on $n_\e$ we first use a Hopf-Cole transformation to deal with bounded values. Next, we prove regularity estimates on variables $x$ and $t$. However, we don't have a priori estimates on variable $y$ due to the nonlocal term $\varrho_\e(t,y)$. To handle this difficulty we first pass to the weak limit, fixing the variable $y$, and next we recover  pointwise and strong convergence from the uniqueness and the structure of the limit.
\\

\noindent
We postpone the statement and proof of a similar result in elliptic case, that is system~\eqref{eq:ntrue} to the end of the paper (Section~\ref{sec:elliptic}). We begin by giving general a priori bounds which hold both for the elliptic and parabolic case, in Section~\ref{sec:prel}. With these at hand, we prove  Theorem~\ref{th:1} in the next section. For the sake of completeness, we recall some H\"older continuity results for parabolic equations in Section~\ref{sec:unifest}. Some conclusions and perspectives are drawn in Section~\ref{sec:conclusion}.

%%%%%%%%%%%%%%%%%%%%
\section{Preliminary estimates}
\label{sec:prel}
%%%%%%%%%%%%%%%%%%%%

Several bounds can be obtained from elementary manipulations of the  equations~\eqref{eq:ntrue} or equations ~\eqref{eq:ne}--\eqref{eq:rho}. Here we make the assumptions~\eqref{as:0}--\eqref{as:2} in the parabolic case and assume~\eqref{as:0} and \eqref{as:2} in the elliptic case. These bounds  are 
\begin{lemma}
The following estimates hold true:
\\[5pt]
(i) $0 \leq c_\e (y,t) \leq c_B$,
\\[5pt]
(ii) $\varrho_\e(y,t)  \leq c_B \dis \max_{0< x<1} \f {r(x)}{d(x)}  -1 = \varrho_M $,
\\[5pt]
(iii)  $c_\e (y,t) \geq c_B \f{\lb}{\lb +\varrho_M} = c_m$,
\\[5pt]
(iv) $ \varrho_\e(y,t)  \geq c_m \dis \min_{0< x<1} \f {r(x)}{d(x)}  -1 = \varrho_m >0$.
\label{lm:est}
\end{lemma}

\proof
We only give the proofs for the parabolic case~\eqref{eq:ne}--\eqref{eq:rho}. The estimates can be proved for the elliptic case~\eqref{eq:ntrue}  following similar arguments.\\ 
We first notice, from \eqref{eq:ne} and \eqref{as:2},  that $n_\e>0$,  in $\R \times \R \times \R^+$. In particular,  $\varrho_\e>0$ in $\R \times \R^+$. Similarly, from \eqref{eq:ce}, \eqref{as:1} and the comparison principle we obtain that $c_\e>0$ in $\R \times \R^+$.
\\ \\
(i)  From  \eqref{eq:ce}, \eqref{as:1}, $\varrho_\e>0$ and the comparison principle, we immediately deduce that $c_\e  \leq c_B$ in $\R \times \R^+$.
\\ \\
(ii) We integrate \eqref{eq:ne} with respect to $x$ and use (i) to obtain
$$
\e \p_t \varrho_\e(y,t) \leq \left( c_B \max_{0< x<1} \f {r(x)}{d(x)}  -  (1+\varrho_\e(y,t)  ) \right) \int d(x) n_\e(x,t) dx .
$$
Using the above inequality together with the definitions \eqref{as:0}, we obtain (ii).
\\
\\
(iii) The third inequality follows directly from (ii), \fer{as:1} and the comparison principle for \fer{eq:ce}.\\
\\
(iv) We integrate \eqref{eq:ne} with respect to $x$, use (iii)  and obtain
$$
\e \p_t \varrho_\e \geq \left( c_m \dis \min_{0< x<1} \f {r(x)}{d(x)} - (1+\varrho_\e ) \right)  \int d(x) n_\e(x,t) dx .
$$
Using the above inequality together with the definitions \eqref{as:0}, we obtain (iv).

%%%%%%%%%%%%%%%%%%%%%%%%%%%%%%%%%%%%%%%%%%%%
\section{The limiting problem}
\label{sec:limit}
%-------------------------------------------
%%%%%%%%%%%%%%%%%%%%%%%%%%%%%%%%%%%%%%%%%%%%

{\em First step. Limits for $\varrho_\e(y,t)$ and $c_\e (y,t)$.} From the uniform bound on $\varrho_\e(y,t) $, we define $\langle \varrho (y,t)\rangle $ as the weak limit 
\beq
\varrho_\e(y,t) \rightharpoonup  \langle \varrho (y,t)\rangle  \quad \text{in }  \;  L^\infty\big(\R \times(0,\infty)\big)\text{-w-}* . 
\label{wlrho}
\eeq 
Then, we may pass to the limit in the equation \eqref{eq:ce} for $c_\e$,  and  for that, we use that
\begin{equation}\label{eq:lc}
c_\e (y,t)  \underset{ \e \to 0 }{\longrightarrow}  c (y,t)  \quad \text{locally uniformly in }  \R \times [0, \infty],
\end{equation}
(see Section~\ref{sec:unifest}  for a proof)  and we find the equation for the limiting nutrient concentration
\begin{equation}\label{eq:c}
\f{\p}{\p t} c-\Delta_y c (y,t) +  [ \,  \langle \varrho (y,t)\rangle +\lb \, ] \; c(y,t) = \lb c_B , \qquad y \in \R .
\end{equation}
Its  solution is $C^{1,\alpha}$ in $y$, with $\alpha\in (0,1)$,  by parabolic regularity (see also Section~\ref{sec:unifest}). 
\bigskip

\noindent {\em Second step. The WKB change of unknown.} Rather than working on $\varrho_\e$ directly, we define  as usual the function
$$
u_\e = \e \ln (n_\e) ,
$$
which satisfies
\begin{equation}\label{eq:ue}
 \partial_t u_\e(y,x,t) =  r(x)  c_\e(t,y)  - d(x) (1+\varrho_\e(y,t)) ,  \qquad y \in \R, \; 0< x <1.
\end{equation}

\noindent
It is standard to derive the bounds 
\begin{equation}\label{est:1}
\begin{cases}
  | \f{\p}{\p t} u_\e(x,y,t) | \leq K(t), \quad  | \f{\p}{\p x} u_\e(x,y,t) | + | \f{\p^3}{\p x^3} u_\e(x,y,t) |+ | \f{\p^3}{\p txx} u_\e(x,y,t) |  \leq K(t),  
\\[8pt]
  \f{\p^2}{\p x^2} u_\e(y,x,t) \leq - a , \quad  u_\e(y,x,t) \leq o(1) .
\end{cases}  
\end{equation}
Indeed, these estimates can be obtained by differentiating \fer{eq:ue} and using \fer{as:2}, \fer{as:3} and Lemma \ref{lm:est}. 
\\

\noindent
For our arguments below, we fix $y$ and  pass to the limit.  Extracting several subsequences which a~priori depend on $y$,
$$
\varrho_\e(y,t) \rightharpoonup \varrho(y,t), \;  L^\infty(0,\infty)\text{-w-}* \qquad \varrho_m \leq \varrho(y,t) \leq \varrho_M, 
$$
$$
u_\e(y,x,t) \longrightarrow u(y,x,t),  \text{ uniformly in } x, \; t\in [0, T], \quad \forall t>0.
$$
Notice that this value $ \varrho(y,t)$ may differ from $\langle \varrho (y,t)\rangle $. Passing to the limit, we find, $y$ by $y$,  that 
\begin{equation} \label{eq:u}
\begin{cases}
u(y,x,t) = u^0(y,x)+ r(x)\int_0^t  c(s,y)ds  - d(x)t +d(x)\int_0^t \varrho(y,s)ds ,  \qquad t \geq 0, \quad 0< x <1, 
\\[10pt]
\dis \max_{0\leq x \leq 1}  u(y,x,t) = 0= u (y, X(y,t), t),
\\[10pt]
u(y,x,t=0) =u^0(y,x).
\end{cases}
\end{equation}
By concavity of $u$ in $y$, the maximum point  $X(y,t)$ is unique.
\\ 

\noindent {\em Third step. $\langle \varrho (y,t)\rangle =  \varrho(y,t)$.} 
Because of the particular structure on the right hand side, we know there is a unique solution $(u(y,t),\int_0^t \varrho (y,s) ds)$ of \eqref{eq:u} for each $y$  (see also \cite{GB.BP:08} for general argument). Therefore the full families $\int_0^t \varrho_\e(y,s) ds$ and $u_\e(y,t)$ converge pointwise and not only subsequences, for each $y$ and $t$. Also by continuous dependence upon the parameter $y$  in the data for \eqref{eq:u}, both $u(y,t)$ and $\int_0^t \varrho (y,s) ds$ also have continuous dependence on $y$.

Consequently $\int_0^t \varrho (y,s) ds = \int_0^t \langle \varrho (y,s) \rangle ds$ and thus $\langle \varrho (y,t)\rangle =  \varrho(y,t)$. However this does not imply strong convergence of $\varrho_\e(y,t) $ in the time variable.
\\

\noindent {\em Fourth  step. The mapping $t \mapsto X(y,t)$ is Lipschitz continuous in $t$, $y$ by $y$.} To prove this, let $X_\e$ be the unique maximum point of $u_\e$ and hence $\p_x u_\e(y,X_\e(y,t),t)=0$. We differentiate this equality with respect to $t$ and find
$$
\p_{xt} u_\e(y,X_\e(y,t),t)+\p_{xx} u_\e(y,X_\e(y,t),t) \f{d}{dt}{X_\e(y,t)}=0.
$$
It follows that
$$
-\p_{xx} u_\e(y,X_\e(y,t),t) \f{d}{dt}{X_\e(y,t)} = r'(x)  c_\e(t,y)  - d'(x) (1+\varrho_\e(y,t)).
$$
As a consequence, $t \mapsto X_\e(y,t)$ is Lipschitz continuous for all $y$.  Therefore, we can pass to the limit as $\e \to 0$, and $X_\e(y,t)$ converges uniformly locally to $X(y,t)$ which is  therefore Lipschitz continuous in $t$ (because this value achieves the maximum of $u(y,t)$ and is unique).  
\\

\noindent
From the estimates on $\partial_{xxx} u$ and $\partial_{txx} u$, we may pass to the strong limit in the term $- \partial_{xx} u_\e\big(y,X_\e(y,t),t \big)$ and obtain 
\beq
\label{X}
\dot X(y,t) =\left(- \partial_{xx} u\big(y,X(y,t),t \big) \right)^{-1} \;  \left( r' \big(X(y,t)\big)  c(t,y)  - d'\big(X(y,t)\big) (1+\varrho(y,t)) \right),
\eeq
Using the arguments in  \cite{GB.BP:08} we can also show, using \fer{eq:u}, that at the Lebesgue points in $t$ of $\varrho(y,t)$ we have
\beq
\label{R0}
 r \big(X(y,t)\big)  c(t,y)  - d \big(X(y,t)\big) (1+\varrho(y,t))  =0. 
\eeq
Since $\varrho\in L^\infty$, we deduce that the above equality holds true for almost every $t$. This implies that $\varrho$ is also H\"{o}lder continuous in $t$. 
\\

\noindent {\em Fifth  step. The mapping $y \mapsto X(y,t)$ is Lipschitz continuous in $y$.} We can use the value of $1+\varrho(y,t)$ given by formula~\eqref{X}, and write the equation~\eqref{R0} for $X$ under the form
$$
\dot X(y,t) =\left(- \partial_{xx} u\big(y,X(y,t),t \big) \right)^{-1} \;  \left( r' \big(X(y,t)\big) - \f{d'\big(X(y,t)\big) r \big(X(y,t)\big)}{d\big(X(y,t)\big)} \right) c(t,y)  .
$$
This is an ordinary differential equation which inherits the regularity of the initial data and coefficients, $c(y,t)$ and $u(y,t)$. Therefore its solution $X(y,t)$ is Lipschitz continuous in $y$. 
\\

\noindent {\em Sixth  step. Strong convergence of $\varrho_\e$.}  With the steps above, the conclusion on the convergence of $n_\e$ is a direct consequence of the analysis of the convergence of $u_\e$. 
\\

\noindent
To conclude, we prove the strong convergence of $\varrho_\e$ following an argument in \cite{AL.SM.BP:10}, equation (9.23). We divide equation \eqref{eq:ne} by $d(x)$ and  integrate. We obtain 
$$
\e \f{\p}{\p t}Ê\int_0^1 \f{n_\e(y,x,t)}{d(x)} dx  =   c_\e(t,y) \int_0^1 \f{r(x) n_\e(y,x,t)}{d(x)}  dx -  \varrho_\e(y,t) (1+\varrho_\e(y,t)).
$$
We pass to the limit and obtain, with $n= \varrho(y,t) \delta \big(x-X(y,t)\big)$, 
$$
 \varrho(y,t)^2 \leq  c(t,y) \int_0^1 \f{r(x) n(y,x,t)}{d(x)}  dx -  \varrho(y,t)=  \varrho(y,t) c(t,y)  \f{r\big(X(y,t) \big) }{d\big(X(y,t) \big)} -  \varrho(y,t).
$$
Comparing with \eqref{R0}, we conclude that this inequality is, in fact, an equality and thus the strong convergence. 
\\

\noindent
The proof of Theorem \ref{th:1} is complete.

%%%%%%%%%%%%%%%%%%%%%%%%%%%%%%%%%%%%%%%%%%%%
\section{Uniform estimates on $c_\e$ (parabolic case)}
\label{sec:unifest}
%-------------------------------------------
%%%%%%%%%%%%%%%%%%%%%%%%%%%%%%%%%%%%%%%%%%%%

In the analysis of the limit of $u_\e$ and $\varrho_\e$, we have used  standard local uniform continuity for $c(y,t)$. We recall the proof for the sake of completeness and show that, locally, $c_\e(y,t)$ is uniformly $1/4$-H\"older continuous in $t$ and  is $1/2$-H\"older continuous in $y$. Better regularity can be obtained using regularizing effects of parabolic equations with the available Lipschitz regularity of $\rho(\cdot,t)$; however, we have chosen to keep a simple complete proof and find a weaker result which is enough for our purpose.
\\
\\ 
{\em First step. Localization method.} We first indicate how to work in $L^2$ after localizing the problem. 
\\

\noindent
Consider a smooth cut-off function $\chi$ with compact support. From equation \eqref{eq:ce} for $c_\e$ (which is uniformly bounded in $L^\infty(\R)$), 
we find
\begin{equation}\label{eq:cchi}
\f{\p}{\p t} [\chi c_\e] -\Delta_y [\chi c_\e (y,t)] + 2 \nabla \chi . \nabla c_\e +   c_\e \Delta \chi = \chi F_\e, \qquad y \in \R , \; t \geq 0,
\end{equation}
with $F_\e = c_B- c_\e (\lb+ \varrho_\e(y,t)) $ which is also  uniformly bounded in $L^\infty$.
\\

\noindent
Therefore, multiplying by $\chi c_\e$ and integrating in $y$, we find after integrations by part 
$$
\f 12  \f d{dt} \int_\R  [\chi c_\e]^2 dy+  \int_\R  |Ê\nabla (\chi c_\e)|^2 dy \leq   \int_\R c_\e^2 \, |\nabla \chi |^2  dy+   \int_\R c_\e F_\e \chi^2 dy \leq C. 
$$
From this estimate, we also control uniformly $\int_0^T \int_B |\nabla c_\e|^2 dy ds$ on each ball $B$ and for each $T\in \R^+$.
\\

\noindent
Next, we multiply equation \eqref {eq:cchi} by $2\Delta_y [\chi c_\e (y,t)]$ and integrate and obtain
$$
\f d{dt} \int_\R  |Ê\nabla( \chi c_\e)|^2 + 2 \int_\R  |\Delta_y (\chi c_\e (y,t))|^2 \leq \int_\R  |\Delta_y (\chi c_\e (y,t))|^2 + R_\e
$$
where $R_\e = \int_\R \big| \chi F_\e - c_\e \Delta \chi -  2 \nabla \chi . \nabla c_\e \big|^2$ is uniformly controlled in $L^1_{\mathrm loc}$  thanks to the  previous estimates on ${\int_0^T} \int_B |\nabla c_\e|^2 dy ds$.
\\

\noindent
As a consequence, for all $0 \leq t \leq T$, 
\begin{equation}\label{eq:para_est}
\int_\R  \big|Ê\nabla [\chi c_\e(y,t) ] \big|^2  dy \leq C_1(T), \qquad \int_0^T \int_\R  \big|\f{\p}{\p t} [\chi c_\e] \big|^2 dy dt \leq C_2(T).
\end{equation}
\\
\\
{\em Second step. H\"older regularity.} We set $v = \chi c_\e$ and prove the following uniform estimate 
%----------------------------
\begin{lemma}
A function $v$ with compact support which satisfies \eqref {eq:para_est} is $1/4$-H\"older continuous in $t$ and  is $1/2$-H\"older continuous in $y$.
\end{lemma}
%------------------------------

\proof The space regularity is obvious since from the first bound (uniform in time) and the Cauchy Schwarz inequality
$$
| v(y_2,t) -v(y_1,t) | \leq \int_{(y_1, y_2)} | \nabla v (y,t)| dy \leq | y_2- y_1 |^{1/2} C_1(T)^{1/2} .
$$

\noindent
Then, we estimate the time increments as follows (with $h= | t_2- t_1 |$) 
$$
| v(y, t_2) -v(y, t_1) | \leq  \int_{(t_1, t_2)} | \f{\p}{\p t}  v (y,t)| dt \leq h^{1/2} \left( \int_0^T | \f{\p}{\p t}  v (y,t)|^2 dt\right)^{1/2}
$$
and thus, being given $y_0$ and $k>0$,  
$$
\int_{|y-y_0| \leq k} | v(y, t_2) -v(y, t_1) | dy \leq  h^{1/2}  \int_{|y-y_0| \leq k}  \left( \int_0^T | \f{\p}{\p t}  v (y,t)|^2 dt\right)^{1/2} dy \leq C_3(T)   h^{1/2}  k^{1/2}.
$$
Finally, we write for all $y \in \R$, 
$$
| v(y_0, t_2) -v(y_0, t_1) | \leq | v(y, t_1) -v(y_0, t_1) | + | v(y_0, t_2) -v(y, t_2) | +  | v(y, t_2) -v(y, t_1) |
$$
and we integrate in $y $ for $|y-y_0|Ê\leq k$ (with $k$ to be chosen later). We find
$$
2 k | v(y_0, t_2) -v(y_0, t_1) | \leq \f{8}{3} C_1(T)^{\f 1 2} k \;  k^{1/2}  + C_3(T) h^{1/2}k^{1/2}.
$$
We take $k=h^{1/2}$ and find the result.

\section{Elliptic coupling}
\label{sec:elliptic}
%-------------------------------------------
%%%%%%%%%%%%%%%%%%%%%%%%%%%%%%%%%%%%%%%%%%%%

In case of elliptic coupling, that is of the system given by equations~\eqref{eq:ntrue}, an additional difficulty occurs because the regularizing effect in time for $c_\e$ and $c$ cannot occur.  Therefore, $c_\e(y,t)$ and $\varrho_\e(y,t)$ have the same regularity in $t$, that is we only handle $L^\infty$ bounds and, consequently,  weak limits. For that reason our result is weaker 
%-------------------------------
\begin{theorem} [Elliptic case]  With assumptions \eqref{as:0},  \eqref{as:1}--\eqref{as:3}, there is $\varrho(y,t)\in  L^\infty\big(\R \times(0,\infty)\big)$, $X(y,t) \in W^{1,\infty}\big(\R\times [0,\infty) \big)$ and $c(y,t) \in L^\infty \big( \R\times (0,\infty)\big) $, such that solutions of \eqref{eq:ntrue} satisfy, as $\e \to 0$, 
\begin{eqnarray}
\varrho_\e \rightharpoonup  \varrho(y,t), \qquad  c_\e \rightharpoonup   c(y,t), \quad  L^\infty\big(\R \times(0,\infty)\big)\text{-w-}*, 
\\
 n_\e(y,x,t) \rightharpoonup   \varrho(y,t) \delta \big(x-X(y,t) \big) , \ \text{ weakly in measures}.
\end{eqnarray}
Moreover, we have, almost everywhere in $t$,  
\begin{equation}\label{eq:cweak}
-\Delta_y c (y,t) +  \langle \varrho \; c \rangle (y,t)  +\lb \; c(y,t) = \lb c_B , \qquad y \in \R, \; t\geq 0 ,
\end{equation}
and
\begin{equation}\label{eq:Xweak}
 r\big(X(y,t)\big) c(y,t) - d\big(X(y,t)\big) \big(1+\varrho(y,t)\big)=0.
 \end{equation}
\label{th:2}
\end{theorem}

\proof We just indicate the modifications to the proof of Theorem~\ref{th:1}. 

\noindent
With the uniform estimates of Lemma~\ref{lm:est}, we can follow the limiting procedure of Section~\ref{sec:limit}.  The limit in \eqref{eq:lc} is just a weak limit because of the time variable, but a.e. in $t$, $c(y,t)$ belongs to $W^{2,\infty} (\R)$ and the equation~\eqref{eq:c} is replaced by equation~\eqref{eq:cweak}.  Then, the analysis of the limit $u(y,t)$ can be performed as in Section~\ref{sec:limit} and both equations~\eqref{X} and \eqref{R0} hold a.e. in $t$.  The $y$-regularity can be  derived for $\varrho$, but not $t$-Lipschitz regularity, because of the lack of time regularity in $c$. However $X$ itself wins one degree of regularity and is indeed Lipschitz continuous. The Sixth step (strong convergence of $\varrho_\e$) also fails. 

%%%%%%%%%%%%%%%%%%%%%%%%%%%%%%%%%%%%%%%%%%%%
\section{Conclusion and perspectives}
\label{sec:conclusion}
%-------------------------------------------

The asymptotic problem we have handled is one of the simplest where both a trait variable $x$ and a space variable $y$ are used. The main difficulty is that the behaviors in these variables are very different because the solution concentrates as a Dirac mass in $x$ and stays bounded in $y$. We do not know of methods adapted to prove compactness in this kind of situations. Indeed, a simple tool would be to prove a priori estimates in Sobolev spaces in the variables $t$ and $y$ for integrals in $x$ as $\rho(y,t)$ here; because of the concentration in the variable  $x$ we cannot expect such regularity except for the Hopf-Cole transform which however gives indirect information. Here we have been able to use uniqueness for the limit in order to recover compactness avoiding strong a priori estimates. This method is limited to the particular situation at hand. We cannot expect it to work in several other situations, for example more general (nonlinear) growth rates under the form $R(y, c, \varrho)$ or dispersion depending on the trait (as in \cite{EB.VC.NM.SM:12,EB.VC:14}).

%%%%%%%%%%%%%%%%%%%%%%%%%%%%%%%%%%%%%%%%%%%%
%
%%%%%%%%%%%%%%%%%%%%%%%%%%%%%%%%%%%
%
%%%%%% BIBLIO %%%%%%%%%%%%%%%%%%%%%%
%
%%%%%%%%%%%%%%%%%%%%%%%%%%%%%%%%%%%%
%\pagestyle{myheadings}

 %%%%%%%%%%%%%%%%%%%%%%%%%%%%%%%%%%%%%%%%%%%%%%%%%%%
 %\bibliography{bibli.bib}
% \bibliographystyle{plain}
 %%%%%%%%%%%%%%%%%%%%%%%%%%%%%%%%%%%%%%%%%%%%%%%%%%%

%\begin{thebibliography}{99}
%\end{thebibliography}

\end{document}